\documentclass{amsart}
\usepackage[english]{babel}
\usepackage[latin1]{inputenc}
\usepackage[dvips,final]{graphics}
\usepackage{amsmath,amsfonts,amssymb,amsthm,amscd,array,
stmaryrd,mathrsfs,mathdots,epigraph}
\usepackage[makeroom]{cancel}
\usepackage{pstricks}
\usepackage[all]{xy}
\usepackage{url}
\usepackage{multirow, blkarray}
\usepackage{booktabs}
\usepackage{caption}
\usepackage{subcaption}

\usepackage{blindtext}
\usepackage{textcomp}
\usepackage[final]{epsfig}
\usepackage{color}
\usepackage{mathabx}

\theoremstyle{plain}
\newtheorem{thm}{Theorem}

\newtheorem{prop}[thm]{Honesty}
\newtheorem*{mthm}{Markov Theorem}

\theoremstyle{definition}

\newcommand{\R}{\mathbb{R}}
\newcommand{\Z}{\mathbb{Z}}

\newcommand{\N}{\mathbb{N}}

\def\a{\alpha}
\def\b{\beta}

\def\e{\varepsilon}
\def\g{\gamma}

\def\vfi{\varphi}

\begin{document}

\title[Shadow sequences of integers]{Shadow sequences of integers,\\
{\small from Fibonacci to Markov and back}}

\author{Valentin Ovsienko}
\address{
Valentin Ovsienko,
Centre National de la Recherche Scientifique,
Laboratoire de Math\'ematiques de Reims, UMR9008 CNRS,
Universit\'e de Reims Champagne-Ardenne,
U.F.R. Sciences Exactes et Naturelles,
Moulin de la Housse - BP 1039,
51687 Reims cedex 2,
France}
\email{valentin.ovsienko@univ-reims.fr}

\maketitle

\thispagestyle{empty}

The Gems and Curiosities column, once called Mathematical Entertainments,
knew its golden years under the watch of David Gale.
For an account, see the amazing book~\cite{Gal1}.
Integer sequences was one of
the main subjects discussed in the column at that time~\cite{Gal}.
The sequences now called Gale-Robinson and Somos sequences first
appeared in the column.
Gale's ``riddles'' about integer sequences strongly influenced combinatorics; 
see for instance the powerful and now classical reference~\cite{FZ}, 
in which Gale's riddles were solved.

As a doomed to fail attempt to maintain the tradition, 
I will discuss a large class of integer sequences given by recurrence relations.
The following general idea looks crazy.
What if another integer sequence follows each integer sequence like a {\it shadow}\/?
I will demonstrate that this is indeed the case,
perhaps not for every integer sequence (this unfortunately I don't know), 
but for many of them.

Does this mean that the number 341962 of known and registered sequences (cf.~\cite{OEIS})
will double? 
The answer is ``no'' for two reasons: shadows of known sequences can also be known,
and one sequence can have several shadows.

\section{Dual numbers}

A pair of integers $A=(a,\a)\in\Z^2$ can be organized as a linear expression
$$
A:=a+\a\e,
$$
where $\e$ is a formal variable.
If~$\e$ is a square root of~$-1$, then~$A$ is a complex number, called a Gaussian integer.
But if~$\e$ satisfies the condition $\e^2=0$,
then~$A$ is called a {\it dual number}. 
Dual numbers were introduced by Clifford in 1873, and have some applications
in geometry and mathematical physics.
In geometry, dual numbers are used to work with the space of oriented lines in $\R^3$, 
a useful device for geometrical optics and computer vision.

Let $(a_n)_{n\in\N}$ be an integer sequence whose every
entry $a_n$ is determined by segment of a fixed length~$k$:
$$
a_{n+k}=R(a_{n+k-1},\ldots,a_{n}),
\qquad\hbox{for all}\quad
n,
$$
where~$R$ is a (generally rational) expression,
 and some initial conditions $(a_1,\ldots,a_k)$.
Assume that a sequence of dual numbers~$(A_n)_{n\in\N}$ satisfies
the same recurrence
$$
A_{n+k}=R(A_{n+k-1},\ldots,A_{n+1}).
$$
The expression in the right-hand-side than has two components:
$R=R_0+R_1\e$,
and $\a_{n+k}$ is then determined by the $\e$-component~$R_1$
and some initial conditions.
Moreover, for some ``mysterious'' reasons, 
the new sequence $(\a_n)_{n\in\N}$ turns out to be integer!

This idea to construct dual integer sequences was suggested in~\cite{Ovs,OT}
and applied to the Gale-Robinson and Somos sequences.
Today I will go further and apply it to several other interesting sequences,
making the method more universal.
The main example is that of the Markov numbers.

\section{The most honest sequence}

The only way to convince whoever that a general method has a chance to work
is to consider examples.
The skeptical reader may interrupt me here and say:
``O.k.! Let us take the `most honest' sequence of all positive integers
$$
a_n=\textcolor{red}{1,2,3,4,5,6,7,\ldots}
$$
registered in~\cite{OEIS} under a (somewhat surprising) number A000027.
What is ithe shadow of A000027?''

Let me see... the sequence satisfies a linear recurrence $a_{n+1}=a_n+1$,
which after substituting to it a sequence of dual numbers $A_n=a_n+\a_n\e$
will produce $\textcolor{blue}{\a,\a,\a,\ldots}$
with arbitrary~$\a$, a sad constant sequence.
But, I promised shadows of sequences, not of poorly chosen recurrence relations!

\begin{prop}
The shadow of A000027 is the sequence A000292 
$$
\a_n=\textcolor{blue}{0,1,4,10,20,35, 56,\ldots}
$$
called the tetrahedral numbers and given explicitly by
$\a_n=\frac{\left(n-1\right) n \left(n+1\right)}{6}$.
\end{prop}

To explain this, notice that, besides the above linear recurrence,
A000027 satisfies another, more interesting non-linear recurrence
\begin{equation}
\label{SimpRec}
a_na_{n+2}=a_{n+1}^2-1.
\end{equation}
Substituting $A_n=a_n+\a_n\e$ instead of~$a_n$ and collecting the $\e$-terms, gives 
the following linear recurrence for $\a_n$:
\begin{equation}
\label{SimpRecBis}
a_n\a_{n+2}=
2a_{n+1}\a_{n+1}-a_{n+2}\a_n.
\end{equation}
More precisely, 
\begin{equation}
\label{PyrRec}
\a_{n+2}=
\frac{2\left(n+1\right)}{n}\a_{n+1}-\frac{\left(n+2\right)}{n}\a_n.
\end{equation}
At first glance, it is not clear that~$\a_n$ stays integer, as~$n$ grows.
But, it is an easy exercice to check that~\eqref{PyrRec} is satisfied for the tetrahedral sequence.
Therefore, choosing the initial conditions $(\a_1,\a_2)=(0,1)$, indeed gives A000292.

Interestingly, the alternative choice $(\a_1,\a_2)=(1,0)$ leads to the sequence
$$
\a_n=\textcolor{blue}{1,0,-3,-9,-19,-34,-55\ldots}
$$
which is (up to a sign) the same sequence A000292 decreased by~$1$
(yet registered as A062748).
An arbitrary solution of~\eqref{PyrRec} is a linear combination of A000292 and A062748.

Another interesting observation is that the sequence of tetrahedral numbers A000292
is actually the {\it convolution} of A000027 with itself.

To end up with the warm-up example of A000027, 
I add that besides~\eqref{SimpRec} it satisfies many other recursions,
for instance, $a_na_{n+3}=a_{n+1}a_{n+2}-2$.
However, the described shadowing procedure does not lead to integer sequences~$\a_n$.
I cannot see any other good candidate for a shadow of A000027 than A000292.

\section{The shadows of Fibonacci and Catalan}

The Fibonacci numbers (see A000045)
$$
F_n=\textcolor{red}{1,\,1,\,2,\,3,\,5,\,8,\,13,\,21,\,34,\,55,\,89,\,144,\,233,\ldots}
$$
satisfy a linear recurrence $F_{n+2}=F_{n+1}+F_n$, 
but, once again, linear recurrence is not interesting for the shadowing purposes.

\begin{prop}
The shadow of the Fibonacci sequence is the sequence A001629 
$$
\vfi_n=\textcolor{blue}{0, 1, 2, 5, 10, 20, 38, 71, 130, 235, 420, 744, 1308,\ldots}
$$
which is the convolution of the Fibonacci sequence with itself.
\end{prop}

Indeed, consider the Cassini identity
$
F_nF_{n+2}=F_{n+1}^2-(-1)^n,
$
called so after the first director of the Paris Observatory,
and substitute into it $\mathcal{F}_n=F_n+\vfi_n\e$.
The recurrence for $\vfi_n$ is then
similar to~\eqref{SimpRecBis}:
$$
\vfi_{n+2}=
\frac{2F_{n+1}\vfi_{n+1}-F_{n+2}\vfi_n}{F_n}.
$$
Once again, the sequence~$(\vfi_n)_{n\in\N}$ is integer for any
choice of the initial conditions.
The initial conditions $(\vfi_1,\vfi_2)=(0,1)$ then lead to A001629,
while any other initial condition provides non-positive sequences.

The sequence of Catalan numbers (see A000108)
$$
C_n=\textcolor{red}{1, 1, 2, 5, 14, 42, 132, 429, 1430, 4862, 16796, 58786,\ldots}
$$
is the same for combinatorics that Fibonacci is for nature,
as they appear in hundreds of combinatorial problems.
The Catalan numbers obey the following recurrence
that was already known to Euler:
$$
C_{n+1}=C_0C_n+C_1C_{n-1}+\cdots+C_nC_0.
$$
Taking the sequence of dual numbers $\mathcal{C}_n=C_n+\g_n\e$
and substituting it into this recurrence, the $\e$-part is determined by
$$
\g_{n+1}=2\left(
C_0\g_n+C_1\g_{n-1}+\cdots+C_n\g_0
\right).
$$
I leave it an exercise to check that,
choosing the initial conditions
$(\g_0,\g_1)=(0,1)$, leads to the sequence A000984
$$
\g_n=\textcolor{blue}{0,1, 2, 6, 20, 70, 252, 924, 3432, 12870, 48620, 184756,\ldots}
$$
of central binomial coefficients $\binom{2n}{n}$.
Choosing the initial conditions
$(\g_0,\g_1)=(1,0)$, leads to A162551
$$
\g_n=\textcolor{blue}{1,0, 2, 8, 30, 112, 420, 1584, 6006, 22880, 87516, 335920,\ldots}
$$
 of double binomials $2\binom{2n}{n-1}$. 
 Both sequences and their linear combinations are good candidates for
 a shadow of the Catalan numbers.

\section{The tree of Markov numbers}

The Markov numbers are triplets of positive integers $(a,b,c)$
which are solutions to the Diophantine equation
\begin{equation}
\label{MarEq}
a^2+b^2+c^2=3abc,
\end{equation}
also called after Markov, who found all the solutions of~\eqref{MarEq}.
Originally related to number theory, the Markov numbers were later found
in geometry and topology~\cite{Rud}, as well as in many other branches
of mathematics; see~\cite{Aig}.
The Markov numbers is still an active area of research 
in combinatorics and mathematical physics.

Let me recall some elements of the Markov theory.
Obviously, $(1,1,1)$ is a solution of~\eqref{MarEq}, and a theorem, 
that Russian mathematicians will always be proud of,
states the following.

\begin{mthm}
Every positive integer solution of~\eqref{MarEq}
can be obtained from $(1,1,1)$ by a sequence of transformations
$(a,b,c)\mapsto(a',b,c)$, where
\begin{equation}
\label{Marmut}
a'=\frac{b^2+c^2}{a},
\end{equation}
and permutations of $a,b,$ and $c$.
\end{mthm}

It is clear that the transformations~\eqref{Marmut} will always produce integers.
Indeed, it follows directly from~\eqref{MarEq}, that it can be rewritten without division
$a'=3bc-a$.
However, for the reasons that I do not explain, 
the form~\eqref{Marmut} is more conceptual.
It is also easy to check that $(a',b,c)$ remains a solution if $(a,b,c)$ is a solution.
The difficult part of the theorem is that every solution is obtained this way.

The Markov numbers can be organized with the help of an infinite binary tree.
The tree is drawn in the plane cutting it into infinitely many regions,
and every region is labeled by a Markov number.
Locally the picture is this
$$
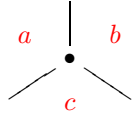
\xymatrix @!0 @R=0.3cm @C=0.3cm
{
&&&\ar@{-}[ddd]\\
\\
&\textcolor{red}{a}&&&&\textcolor{red}{b}\\
&&&\bullet\ar@{-}[llldd]\ar@{-}[rrrdd]\\
\\
&&&\textcolor{red}{c}&&&
}
$$
and the transformations~\eqref{Marmut} correspond to the following branchings:
$$
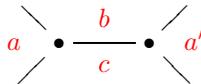
\xymatrix @!0 @R=0.3cm @C=0.3cm
{
\ar@{-}[rrdd]&&&&&&&&\\
&&&&\textcolor{red}{b}\\
\textcolor{red}{a}&&\bullet\ar@{-}[lldd]\ar@{-}[rrrr]&&&&
\bullet\ar@{-}[rruu]\ar@{-}[rrdd]&&\textcolor{red}{a'}\\
&&&&\textcolor{red}{c}\\
&&&&&&&&
}
$$
The tree of Markov numbers grows like this
$$
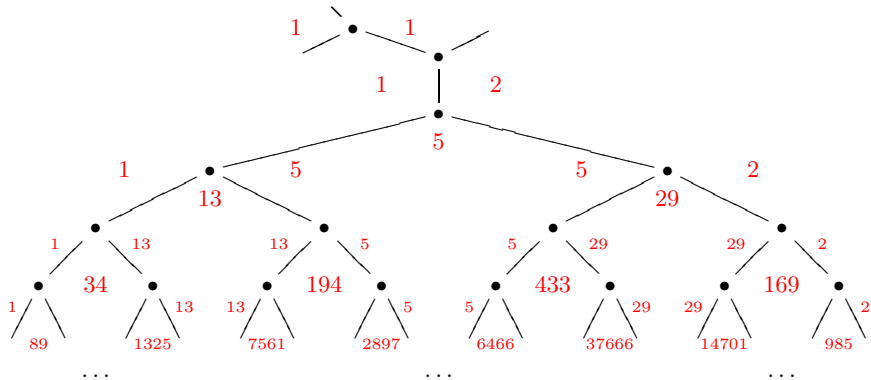
\begin{small}
\xymatrix @!0 @R=0.38cm @C=0.38cm
{&&&&&&&&&&&&\\
&&&&&&&&&&&{\textcolor{red}{1}}&&\bullet\ar@{-}[lld]\ar@{-}[lu]&&{\textcolor{red}{1}}&&&\\
&&&&&&&&&&&&&&&&\bullet\ar@{-}[dd]\ar@{-}[rru]\ar@{-}[lllu]\\
&&&&&&&&&&&&&&{\textcolor{red}{1}}&&&&{\textcolor{red}{2}}\\
&&&&&&&&&&&&&&&&\bullet\ar@{-}[lllllllldd]\ar@{-}[rrrrrrrrdd]&&&&&&&&\\
&&&&&&&&&&&&&&&&{\textcolor{red}{5}}\\
&&&&&{\textcolor{red}{1}}&&&\bullet\ar@{-}[lllldd]\ar@{-}[rrrrdd]
&&&{\textcolor{red}{5}}&&&&&&&&&&{\textcolor{red}{5}}&&&\bullet\ar@{-}[lllldd]\ar@{-}[rrrrdd]&&&{\textcolor{red}{2}}\\
&&&&&&&&{\textcolor{red}{13}}
&&&&&&&&&&&&&&&&{\textcolor{red}{29}}\\
&&&&\bullet\ar@{-}[lldd]_{\textcolor{red}{1}}\ar@{-}[rrdd]^{\textcolor{red}{13}}
&&&&&&&&\bullet\ar@{-}[lldd]_{\textcolor{red}{13}}\ar@{-}[rrdd]^{\textcolor{red}{5}}
&&&&&&&&\bullet\ar@{-}[lldd]_{\textcolor{red}{5}}\ar@{-}[rrdd]^{\textcolor{red}{29}}
&&&&&&&&\bullet\ar@{-}[lldd]_{\textcolor{red}{29}}\ar@{-}[rrdd]^{\textcolor{red}{2}}\\
&&&&&&&&&&&&&&&&&&&&&&&&&&&&\\
&&\bullet\ar@{-}[ldd]_{\textcolor{red}{1}}\ar@{-}[rdd]
&&{\textcolor{red}{34}}&&\bullet\ar@{-}[ldd]\ar@{-}[rdd]^{\textcolor{red}{13}}
&&&&\bullet\ar@{-}[ldd]_{\textcolor{red}{13}}\ar@{-}[rdd]
&&{\textcolor{red}{194}}&&\bullet\ar@{-}[ldd]\ar@{-}[rdd]^{\textcolor{red}{5}}
&&&&\bullet\ar@{-}[ldd]_{\textcolor{red}{5}}\ar@{-}[rdd]
&&{\textcolor{red}{433}}&&\bullet\ar@{-}[ldd]\ar@{-}[rdd]^{\textcolor{red}{29}}
&&&&\bullet\ar@{-}[ldd]_{\textcolor{red}{29}}\ar@{-}[rdd]
&&{\textcolor{red}{169}}&&\bullet\ar@{-}[ldd]\ar@{-}[rdd]^{\textcolor{red}{2}}\\
&&&&&&&&&&&&&&&&&&&&&&
&&&&&&&&&&&\\
&&{\textcolor{red}{\scriptstyle89}}
&&&&{\textcolor{red}{\scriptstyle1325}}
&&&&{\textcolor{red}{\scriptstyle7561}}
&&&&{\textcolor{red}{\scriptstyle2897}}
&&&&{\textcolor{red}{\scriptstyle6466}}
&&&&{\textcolor{red}{\scriptstyle37666}}
&&&&{\textcolor{red}{\scriptstyle14701}}
&&&&{\textcolor{red}{\scriptstyle985}}&&\\
&&&&\ldots&&&&&&&&&&&&\ldots&&&&&&&&&&&&\ldots
}
\end{small}$$

The simplest branches of the Markov tree are those bounded by~$1$ and~$2$,
the border branches in the above picture.
The corresponding triplets of Markov numbers contain the odd Fibonacci numbers 
$(1,F_{2k-1},F_{2k+1})$; see A001519, and
the odd Pell numbers $(2,P_{2k-1},P_{2k+1})$; see A001653.
Various subsequences and arrangements of the Markov numbers
are registered in the OEIS as dozens of entries; see A002559 and related sequences.

\section{The shadow of Andrey Andreyevich Markov}
My ultimate goal is everything but to cast a shadow over the Markov numbers,
but it is tempting to apply the general shadowing construction!..
It goes very simply:
choose the initial triplet of dual numbers $(A_1,B_1,C_1)$ with
\begin{equation}
\label{InitD}
A_1=1+\a_1\e,\qquad
B_1= 1+\b_1\e,\qquad
C_1=1+\g_1\e,
\end{equation}
where $\a_1,\b_1,\g_1$ are arbitrary integers,
and then apply all possible sequences of the transformations $(A,B,C)\mapsto(A',B,C)$, with
$$
A'=\frac{B^2+C^2}{A},
$$
mixed with permutations of $A,B$ and~$C$.
Collecting the terms with and without~$\e$, one obtains
\begin{equation}
\label{MarmutD}
a'=\frac{b^2+c^2}{a},
\qquad\qquad
\a'=\frac{2b\,\b+2c\,\g-a'\a}{a}.
\end{equation}
This time, it is not at all obvious why $\a'$ will remain integer.
\begin{prop}
For an arbitrary choice of the initial conditions~\eqref{InitD},
the transformations~\eqref{MarmutD}, mixed with permutations,
produce integer sequences $(\a_n)_{n\in\N},(\b_n)_{n\in\N}$ and~$(\g_n)_{n\in\N}$.
\end{prop}

This integrality persists thanks to a ``miracle'', called the Laurent phenomenon~\cite{FZ},
the same miracle that guarantees integrality of the Somos and Gale-Robinson sequences.
The complete proof is too technical to be reproduced here, a general statement can be found in~\cite{OT}.

There is one choice of the initial values $(\a_1,\b_1,\g_1)$ 
which seems to be natural and interesting:
\begin{equation}
\label{InitGooD}
(\a_1,\b_1,\g_1)=(0,1,1).
\end{equation}
Let me explain this.
First, it is natural to take it symmetric in~$b$ and~$c$, that is,
to assume $\b_1=\g_1$.
Choosing $(\a_1,\b_1,\g_1)=(1,1,1)$ would produce the same Markov numbers,
in the sense that the $\e$-part would coincide with the classical one:
$$
\a_n=a_n,
\qquad
\b_n=b_n,
\qquad
\g_n=c_n
$$ 
which is not very interesting.
Choosing $(\a_1,\b_1,\g_1)=(1,0,0)$ would produce negative numbers.
The choice~\eqref{InitGooD} is the only remaining.

Here is the tree of Markov numbers together with its shadow
$$
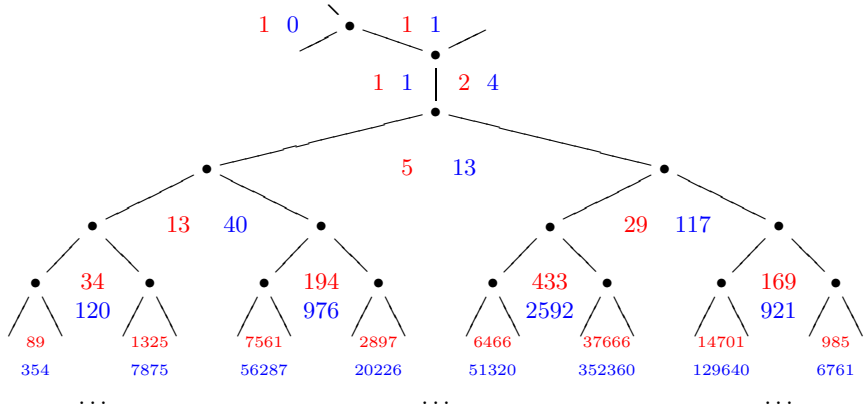
\begin{small}
\xymatrix @!0 @R=0.38cm @C=0.38cm
{&&&&&&&&&&&&\\
&&&&&&&&&&{\textcolor{red}{1}}&{\textcolor{blue}{0}}
&&\bullet\ar@{-}[lld]\ar@{-}[lu]&&{\textcolor{red}{1}}&{\textcolor{blue}{1}}&&\\
&&&&&&&&&&&&&&&&\bullet\ar@{-}[dd]\ar@{-}[rru]\ar@{-}[lllu]\\
&&&&&&&&&&&&&&{\textcolor{red}{1}}&{\textcolor{blue}{1}}&&{\textcolor{red}{2}}&{\textcolor{blue}{4}}\\
&&&&&&&&&&&&&&&&\bullet\ar@{-}[lllllllldd]\ar@{-}[rrrrrrrrdd]&&&&&&&&\\
&&&&&&&&&&&&&&&&\\
&&&&&&&&\bullet\ar@{-}[lllldd]\ar@{-}[rrrrdd]
&&&&&&&{\textcolor{red}{5}}&&{\textcolor{blue}{13}}
&&&&&&&\bullet\ar@{-}[lllldd]\ar@{-}[rrrrdd]&&&\\
&&&&&&&&
&&&&&&&&&&&&&&&&\\
&&&&\bullet\ar@{-}[lldd]\ar@{-}[rrdd]
&&&{\textcolor{red}{13}}&&{\textcolor{blue}{40}}&&&\bullet\ar@{-}[lldd]\ar@{-}[rrdd]
&&&&&&&&\bullet\ar@{-}[lldd]\ar@{-}[rrdd]
&&&{\textcolor{red}{29}}&&{\textcolor{blue}{117}}&&&\bullet\ar@{-}[lldd]\ar@{-}[rrdd]\\
&&&&&&&&&&&&&&&&&&&&&&&&&&&&\\
&&\bullet\ar@{-}[ldd]\ar@{-}[rdd]
&&{\textcolor{red}{34}}&&\bullet\ar@{-}[ldd]\ar@{-}[rdd]
&&&&\bullet\ar@{-}[ldd]\ar@{-}[rdd]
&&{\textcolor{red}{194}}&&\bullet\ar@{-}[ldd]\ar@{-}[rdd]
&&&&\bullet\ar@{-}[ldd]\ar@{-}[rdd]
&&{\textcolor{red}{433}}&&\bullet\ar@{-}[ldd]\ar@{-}[rdd]
&&&&\bullet\ar@{-}[ldd]\ar@{-}[rdd]
&&{\textcolor{red}{169}}&&\bullet\ar@{-}[ldd]\ar@{-}[rdd]\\
&&&&{\textcolor{blue}{120}}&&&&&&&&{\textcolor{blue}{976}}&&&&&&&&{\textcolor{blue}{2592}}&&
&&&&&&{\textcolor{blue}{921}}&&&&&\\
&&{\textcolor{red}{\scriptstyle89}}
&&&&{\textcolor{red}{\scriptstyle1325}}
&&&&{\textcolor{red}{\scriptstyle7561}}
&&&&{\textcolor{red}{\scriptstyle2897}}
&&&&{\textcolor{red}{\scriptstyle6466}}
&&&&{\textcolor{red}{\scriptstyle37666}}
&&&&{\textcolor{red}{\scriptstyle14701}}
&&&&{\textcolor{red}{\scriptstyle985}}&&\\
&&{\textcolor{blue}{\scriptstyle354}}
&&&&{\textcolor{blue}{\scriptstyle7875}}
&&&&{\textcolor{blue}{\scriptstyle56287}}
&&&&{\textcolor{blue}{\scriptstyle20226}}
&&&&{\textcolor{blue}{\scriptstyle51320}}
&&&&{\textcolor{blue}{\scriptstyle352360}}
&&&&{\textcolor{blue}{\scriptstyle129640}}
&&&&{\textcolor{blue}{\scriptstyle6761}}
\\
&&&&\ldots&&&&&&&&&&&&\ldots&&&&&&&&&&&&\ldots
}
\end{small}$$
Some observations can be made.
The sequence 
$$
{\textcolor{blue}{1,4,13,40,120,354,1031,2972, 8495,\ldots}}
$$
appearing as a compagnon of the odd Fibonacci branch,
turns out to be known.
This is the delightful A238846 which is the convolution of two
bisections of the Fibonacci sequence, $F_{2n+1}$ and $F_{2n}$.
Other subsequences appear to be new.

What is the role of the shadow Markov tree? 
Does it mean something?
I would be glad if I could solve this riddle,
but at this stage, I can only say:

\hfill{
Cottleston, Cottleston, Cottleston Pie,}

\hfill{
A fish can't whistle and neither can I.}

\hfill{
Ask me a riddle and I reply:}

\hfill{
``Cottleston, Cottleston, Cottleston Pie.}\footnote{
from AA. Milne, Winnie the Pooh.}''

\end{document}